\documentclass[12pt,a4paper]{article}
\usepackage{latexsym,amsfonts,graphicx}

\newcommand{\G}{\Gamma}
\newcommand{\tg}{\Gamma}
\newcommand{\htr}{\Delta}

\newcommand{\F}{{\cal F}}
\renewcommand{\L}{{\cal L}}

\newcommand{\C}{\mathbb{C}}

\begin{document}

\title{Note on the Schwarz Triangle functions}
\author{Mark Harmer\\
email: m.s.harmer@massey.ac.nz\\
Institute of Information and Mathematical Sciences\\
Massey University Albany\\
New Zealand}

\maketitle

\begin{abstract}
We shown the rationality of the Taylor coefficients of the inverse of the Schwarz triangle functions for a triangle group about any vertex of the fundamental domain.
\end{abstract}

\section{Introduction}
In the paper \cite{Leh2} of the same title Lehner shows that the (inverse of the) Schwarz triangle functions for the Hecke triangle groups have rational Fourier coefficients about the cusp along with an asymptotic formula for the coefficients. Here we generalise this result to show that the Taylor coefficients of any Schwarz triangle function are rational about any vertex of the fundamental domain (in appropriate coordinates). \\
We consider the triangle groups $\tg$; ie. groups generated by pairs of reflections across the sides of the hyperbolic triangle $\htr$ with angles $\pi\alpha = \frac{\pi}{m}$, $\pi\beta = \frac{\pi}{n}$ and $\pi\gamma = \frac{\pi}{p}$ at the vertices for positive integers $m$, $n$ and $p$ with $\frac{1}{m} + \frac{1}{n} + \frac{1}{p} <1$ \cite{Bea}. A fundamental domain $\F$ of $\tg$ consists of the union of a copy of $\htr$ with its reflection through any one of its sides. \\
The Riemann mapping theorem tells us that it is possible to find a conformal map $\varphi (z)=w$ from the upper half plane to $\htr$---the classical Schwarz triangle functions \cite{Neh}. The inverse function $\psi(w)=z$ is a single valued meromorphic function which is automorphic with respect to the group, ie. $\psi(w)=\psi(g w)$ for all $g\in\tg$. Here we prove the simple result that $\psi(w)$ has rational Taylor coefficients when expanded about any of the vertices of the fundamental domain in appropriate coordinates. \\
\section{Schwarz-Christoffel Mappings to Regions enclosed by Circular Arcs}
\begin{figure}[ht]\hspace*{-45mm}
\includegraphics{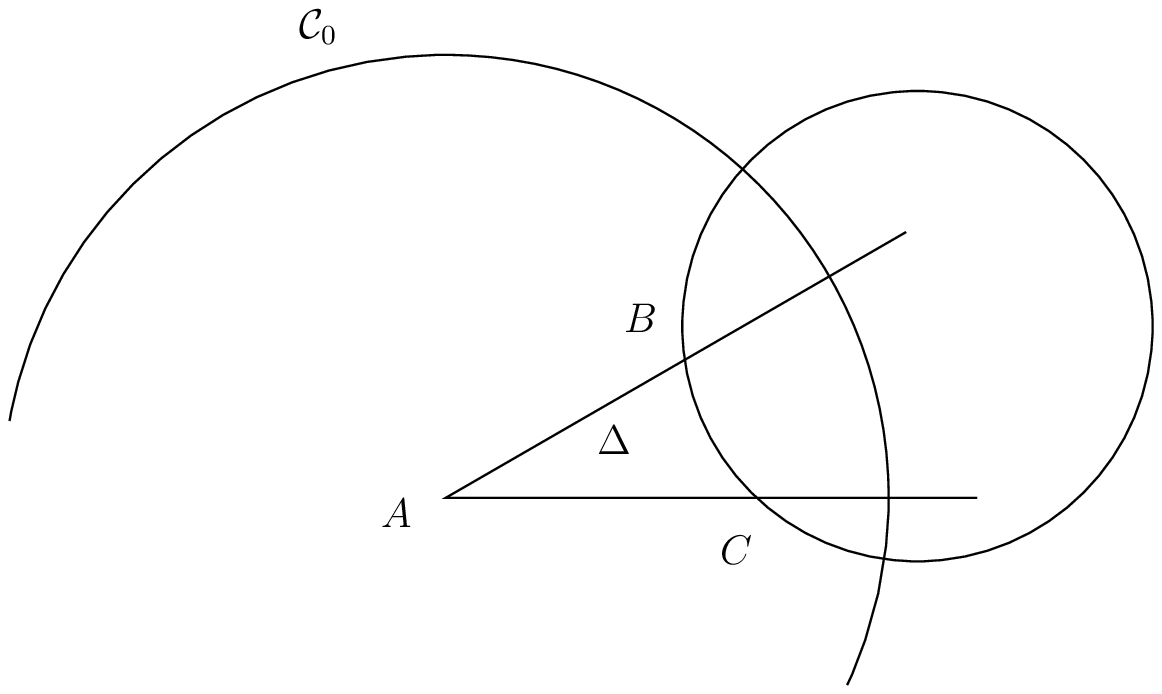}
\vspace*{-198mm}
\caption{The image $\htr$ of $\varphi$.}\label{orthoc}
\end{figure}

A detailed discussion of the theory of Schwarz-Christoffel mappings may be found in
\cite{Neh}. To summarise, given a region $\htr\subset\C$ enclosed by three circular arcs with angles $\left\{ \pi \alpha, \pi\beta, \pi\gamma \right\}$ at the respective vertices $\left\{ A, B , C \right\}$ we can find an holomorphic function
$$
\varphi (z) = \frac{z^{1-c} F( a',b';c';z )}{F( a, b; c; z )}
$$
which maps the upper half plane $\C_+=\{ z;\; \Im z>0\}$ onto the interior of $\htr$ and is continuous to the boundary. Here the numerator and denominator are linearly independent solutions of the hypergeometric equation
$$
z(1-z) F_{zz} + [c-(a+b+1)z] F_{z} - ab F = 0 
$$
with constants $a$, $b$, $c$ given in terms of the angles by 
\begin{eqnarray*}
& & a = {\textstyle \frac{1}{2}} (1-\alpha+\beta-\gamma) \label{prms1} \\
& & b = {\textstyle \frac{1}{2}} (1-\alpha-\beta-\gamma) \label{prms2} \\
& & c = 1-\alpha \, . \label{prms3}
\end{eqnarray*}
The function $F$ is the hypergeometric function and the arguments in the numerator are given by
\begin{eqnarray*}
& & a' = a-c+1 = {\textstyle \frac{1}{2}} (1+\alpha+\beta-\gamma) \label{aprms1} \\
& & b' = b-c+1 = {\textstyle \frac{1}{2}} (1+\alpha-\beta-\gamma) \label{aprms2} \\
& & c' = 2-c = 1+\alpha \, . \label{aprms3} 
\end{eqnarray*}
The vertices $A$, $B$ and $C$ of $\htr$ are the images of points $A=\varphi(0)$, $B=\varphi(\infty)$
and $C=\varphi(1)$ on the real axis. Then it is easy to see that $A=\varphi(0)=0$ lies at the origin, $C$ lies on the positive real axis and the arcs $\overline{AB}$, $\overline{AC}$ are straight lines, see \cite{Neh} for details. In fact $\htr$ is a hyperbolic triangle with the vertex $A$ at the origin of a {\em scaled} Poincar\'{e} disc model of the hyperbolic
plane, see figure \ref{orthoc} (here the orthogonal circle ${\cal C}_0$ does not have radius one). \\
The normalisation $\hat{\varphi} (z) = \nu \varphi (z)$ so that ${\cal C}_0$ has radius one and $\htr$ is a hyperbolic triangle in the standard Poincar\'{e} disc is given by 
$$
\nu=\sqrt{\frac{\cos(\pi\alpha+\pi\beta) + \cos(\pi\gamma)}
{\cos(\pi\alpha-\pi\beta) + \cos(\pi\gamma)} \cdot \frac{\cos(\pi\alpha-\pi\beta-\pi\gamma ) + 1}{\cos(\pi\alpha+\pi\beta+\pi\gamma ) + 1} }\cdot \frac{\G (a')\G (b')}{\G (c')} \cdot \frac{\G (c)}{\G (a)\G (b)}  \, . 
$$
For details of this messy calculations see \cite{Har7}. It is clear from this that $\hat{\varphi} (z)=G(\alpha;z)/G(-\alpha;z)$ where $G(\alpha;z) = \sqrt{\nu z^{\alpha}} F( a',b';c';z )$. It is interesting to note that $G$ is the solution of an eigenvalue equation $\L_{\beta,\gamma} G = \alpha^2 G$ with eigenvalue $\alpha^2$ where $\L_{\beta,\gamma}$ is an operator proportional to the hypergeometric operator and independent of $\alpha$.
\section{Rationality of the Taylor Coefficients}
The function $\varphi (z)=w$ is a Schwarzian triangle function. We may analytically continue $\varphi$ through the real axis, avoiding the singular points $\{ 0,1,\infty\}$ which correspond to vertices of $\htr$, to get a map from the cut Riemann sphere to the union of $\htr$ and $\htr$ reflected through one of its sides---together forming a fundamental domain of $\tg$. Continuing in this way $\varphi$ is an infinitely many-valued function on the Riemann sphere minus $\{ 0,1,\infty\}$. On the other hand the inverse function $\psi(w)=z$ is a single valued meromorphic automorphic form: there is a single pole at $B$, and for a suitable definition of the group action it is automorphic, $\psi(g\cdot w) = \psi(w)$ for all $g\in\tg$ (\cite{Neh}, pg. 308). \\
We now show that, in suitable coordinates, the Taylor coefficients of $\psi (w)$ at each of the vertices of the fundamental domain $\F$ are rational. In our argument we will assume that each of $\{m,n,p\}$ are finite however this is not essential---the case of a cusp having been dealt with in \cite{Leh2}. \\
We begin by considering the Taylor series of $\psi (w)$ about the order $m$ vertex $A = 0$; due to automorphy it has the form 
$$
z=\psi (w) = w^{m} + c_2 w^{2m} + c_3 w^{3m} + \cdots \, .
$$
Conversely the Schwarz triangle function has an order $m$ branch point at the origin
\begin{eqnarray*}
w=\varphi (z) & = & z^{\frac{1}{m}} \left( 1+ d_1 z + d_2 z^2 + \cdots \right) \\
& = & z^{\frac{1}{m}} \frac{F (a',b';c';z)}{F (a,b;c;z)} \, .
\end{eqnarray*}
Now the $a$, $a'$ etc. are rational and expanding $F$ in a hypergeometric series about the origin we see that as a result all the $d_i$ will be rational. This in turn implies that in
\begin{eqnarray*}
w^{m} & = & z \left( 1 + d'_1 z + d'_2 z^2 + \cdots \right) \\
w^{2m} & = & z^2 \left( 1 + d''_1 z + d''_2 z^2 + \cdots \right) \\
\vdots &  & \\
w^{km} & = & z^k \left( 1 + d^{(k)}_1 z + d^{(k)}_2 z^2 + \cdots \right) 
\end{eqnarray*}
the $d^ {(k)}_i$ are also rational. Substituting for the powers of $w^{m}$ in $\psi (w)$ we solve for the $c_i$ and see that they too are rational.\\
For the other vertices we define suitable coordinates about the vertex using a linear fractional transformation. Let us bracket $\psi (w)$ with linear transformations $\psi' (w) = T \circ \psi \left( S^{-1} \circ w\right)$ where 
$$
S\circ w = e^{it} \frac{w-d}{1-\tilde{d}w}
$$
maps the hyperbolic triangle in figure \ref{orthoc} so that the vertex $B$ (in turn $C$) is placed at the origin of the (scaled) Poincar\'{e} disc model and the triangle is rotated so that one side is on the real axis and one vertex in the first quadrant. Here $|d|<\nu^{-1}$ and $\tilde{d} = \bar{d}\nu^{2}$ as we use the scaled Poincar\'{e} disc. Then we put
$$
T\circ z = \frac{-1}{z-1}\; \; \; \left(  \mbox{in turn}\; \frac{z-1}{z} \right)
$$
which is the shift $\{ 0,1,\infty\}\to\{ 1,\infty ,0\}$, (in turn $\{ 0,1,\infty\}\to\{ \infty ,0,1\}$). We claim that $\psi' (w)$ is the inverse of a Schwarz triangle function constructed as above with the angles shifted $\{ \alpha, \beta, \gamma \}\to\{ \beta, \gamma, \alpha \}$ (or $\{ \alpha, \beta, \gamma \}\to\{ \gamma, \alpha, \beta \}$). Indeed it is easy to see that $\psi' (w)$ has the same action mapping the upper half plane to the appropriate hyperbolic triangle with the distinguished points $\{ 0,1,\infty\}$ going to the vertices. The uniqueness of this conformal map follows from the Riemann mapping theorem and the fact that we are mapping distinguished points to the vertices---so that $\psi' (w)$ must correspond to the inverse of the constructed Schwarz triangle function. Consequently, we apply the above argument to get the rationality of the Taylor coefficients of 
$$
\psi \left( S^{-1} \circ w\right)
$$ 
about $B$ (and likewise $C$)---it is clear that $T$ will not effect rationality. In these coordinates the Taylor coefficients of the automorphic function $\psi (w)$ are rational about each of the vertices of the fundamental domain.



\end{document}